\magnification\magstep1
\newread\AUX\immediate\openin\AUX=\jobname.aux
\def\ref#1{\expandafter\edef\csname#1\endcsname}
\ifeof\AUX\immediate\write16{\jobname.aux gibt es nicht!}\else
\input \jobname.aux
\fi\immediate\closein\AUX
\def\today{\number\day.~\ifcase\month\or
  Januar\or Februar\or M\"arz\or April\or Mai\or Juni\or
  Juli\or August\or September\or Oktober\or November\or Dezember\fi
  \space\number\year}
\font\sevenex=cmex7
\scriptfont3=\sevenex
\font\fiveex=cmex10 scaled 500
\scriptscriptfont3=\fiveex

\def\P{{\bf P}}

\def\phi{\varphi}
\def\epsilon{\varepsilon}
\def\theta{\vartheta}

\def\uauf{\lower1.7pt\hbox to 3pt{%
\vbox{\offinterlineskip
\hbox{\vbox to 8.5pt{\leaders\vrule width0.2pt\vfill}%
\kern-.3pt\hbox{\lams\char"76}\kern-0.3pt%
$\raise1pt\hbox{\lams\char"76}$}}\hfil}}
\def\cite#1{\expandafter\ifx\csname#1\endcsname\relax
{\bf?}\immediate\write16{#1 ist nicht definiert!}\else\csname#1\endcsname\fi}
\def\expandwrite#1#2{\edef\next{\write#1{#2}}\next}
\def\neverexpand{\noexpand\noexpand\noexpand}
\def\strip#1\ {}
\def\ncite#1{\expandafter\ifx\csname#1\endcsname\relax
{\bf?}\immediate\write16{#1 ist nicht definiert!}\else
\expandafter\expandafter\expandafter\strip\csname#1\endcsname\fi}
\newwrite\AUX
\immediate\openout\AUX=\jobname.aux
\newcount\Abschnitt\Abschnitt0
\def\beginsection#1. #2 \par{\advance\Abschnitt1%
\vskip0pt plus.10\vsize\penalty-250
\vskip0pt plus-.10\vsize\bigskip\vskip\parskip
\edef\TEST{\number\Abschnitt}
\expandafter\ifx\csname#1\endcsname\TEST\relax\else
\immediate\write16{#1 hat sich geaendert!}\fi
\expandwrite\AUX{\neverexpand\ref{#1}{\TEST}}
\leftline{\bf\number\Abschnitt. \ignorespaces#2}%
\nobreak\smallskip\noindent\SATZ1}
\def\Proof:{\par\noindent{\it Proof:}}
\def\Remark:{\ifdim\lastskip<\medskipamount\removelastskip\medskip\fi
\noindent{\bf Remark:}}
\def\Remarks:{\ifdim\lastskip<\medskipamount\removelastskip\medskip\fi
\noindent{\bf Remarks:}}
\def\Definition:{\ifdim\lastskip<\medskipamount\removelastskip\medskip\fi
\noindent{\bf Definition:}}
\def\Example:{\ifdim\lastskip<\medskipamount\removelastskip\medskip\fi
\noindent{\bf Example:}}
\newcount\SATZ\SATZ1
\def\proclaim #1. #2\par{\ifdim\lastskip<\medskipamount\removelastskip
\medskip\fi
\noindent{\bf#1.\ }{\it#2}\Par
\ifdim\lastskip<\medskipamount\removelastskip\goodbreak\medskip\fi}
\def\Aussage#1{%
\expandafter\def\csname#1\endcsname##1.{\ifx?##1?\relax\else
\edef\TEST{#1\penalty10000\ \number\Abschnitt.\number\SATZ}
\expandafter\ifx\csname##1\endcsname\TEST\relax\else
\immediate\write16{##1 hat sich geaendert!}\fi
\expandwrite\AUX{\neverexpand\ref{##1}{\TEST}}\fi
\proclaim {\number\Abschnitt.\number\SATZ. #1\global\advance\SATZ1}.}}
\Aussage{Theorem}
\Aussage{Proposition}
\Aussage{Corollary}
\Aussage{Lemma}
\font\la=lasy10
\def\strich{\hbox{$\vcenter{\hbox
to 1pt{\leaders\hrule height -0,2pt depth 0,6pt\hfil}}$}}
\def\dashedrightarrow{\hbox{%
\hbox to 0,5cm{\leaders\hbox to 2pt{\hfil\strich\hfil}\hfil}%
\kern-2pt\hbox{\la\char\string"29}}}

\def\Bindestrich{\penalty10000-\hskip0pt}
\let\_=\Bindestrich
\def\.{{\sfcode`.=1000.}}

\def\Par{\par}
\def\:={\mathrel{\raise0,9pt\hbox{.}\kern-2,77779pt
\raise3pt\hbox{.}\kern-2,5pt=}}
\def\=:{\mathrel{=\kern-2,5pt\raise0,9pt\hbox{.}\kern-2,77779pt
\raise3pt\hbox{.}}}

\def\Ugleich{\hbox{$\cup$\kern.5pt\vrule depth -0.5pt}}
\def\|#1|{\mathop{\rm#1}\nolimits}
\def\<{\langle}
\def\>{\rangle}
\let\Times=\times
\def\times{\mathop{\Times}}
\let\Otimes=\otimes
\def\otimes{\mathop{\Otimes}}
\catcode`\@=11
\def\hex#1{\ifcase#1 0\or1\or2\or3\or4\or5\or6\or7\or8\or9\or A\or B\or
C\or D\or E\or F\else\message{Warnung: Setze hex#1=0}0\fi}
\def\fontdef#1:#2,#3,#4.{%
\alloc@8\fam\chardef\sixt@@n\FAM
\ifx!#2!\else\expandafter\font\csname text#1\endcsname=#2
\textfont\the\FAM=\csname text#1\endcsname\fi
\ifx!#3!\else\expandafter\font\csname script#1\endcsname=#3
\scriptfont\the\FAM=\csname script#1\endcsname\fi
\ifx!#4!\else\expandafter\font\csname scriptscript#1\endcsname=#4
\scriptscriptfont\the\FAM=\csname scriptscript#1\endcsname\fi
\expandafter\edef\csname #1\endcsname{\fam\the\FAM\csname text#1\endcsname}
\expandafter\edef\csname hex#1fam\endcsname{\hex\FAM}}
\catcode`\@=12 

\fontdef Ss:cmss10,,.
\fontdef Fr:eufm10,eufm7,eufm5.


%
\newread\AUXX
\immediate\openin\AUXX=msxym.tex
\ifeof\AUXX
\fontdef bbb:msbm10,msbm7,msbm5.
\fontdef mbf:cmmib10,cmmib7,.
\else
\fontdef bbb:msym10,msym7,msym5.
\fontdef mbf:cmmib10,cmmib10 scaled 700,.
\fi
\immediate\closein\AUXX
\def\CC{{\bbb C}}
\def\FF{{\bbb F}}

\def\QQ{{\bbb Q}}

\def\ZZ{{\bbb Z}}
\mathchardef\leer=\string"0\hexbbbfam3F
\mathchardef\subsetneq=\string"3\hexbbbfam24
\mathchardef\semidir=\string"2\hexbbbfam6E
\mathchardef\dirsemi=\string"2\hexbbbfam6F
\let\OL=\overline
\def\overline#1{{\hskip1pt\OL{\hskip-1pt#1\hskip-1pt}\hskip1pt}}


%
\abovedisplayskip 9.0pt plus 3.0pt minus 3.0pt
\belowdisplayskip 9.0pt plus 3.0pt minus 3.0pt
\newdimen\Grenze\Grenze2\parindent\advance\Grenze1em
\newdimen\Breite
\newbox\DpBox
\def\NewDisplay#1$${\Breite\hsize\advance\Breite-\hangindent
\setbox\DpBox=\hbox{\hskip2\parindent$\displaystyle{#1}$}%
\ifnum\predisplaysize<\Grenze\abovedisplayskip\abovedisplayshortskip
\belowdisplayskip\belowdisplayshortskip\fi
\global\futurelet\nexttok\WEITER}
\def\WEITER{\ifx\nexttok\qed\expandafter\leftQEDdisplay
\else\leftdisplay\fi}
\def\leftdisplay{\hskip-\hangindent\leftline{\box\DpBox}$$}
\def\leftQEDdisplay{\hskip-\hangindent
\line{\copy\DpBox\hfill\lower\dp\DpBox\copy\QEDbox}%
\belowdisplayskip0pt$$\bigskip\let\nexttok=}
\everydisplay{\NewDisplay}
\newbox\QEDbox
\newbox\nichts\setbox\nichts=\vbox{}\wd\nichts=2mm\ht\nichts=2mm
\setbox\QEDbox=\hbox{\vrule\vbox{\hrule\copy\nichts\hrule}\vrule}
\def\qed{\leavevmode\unskip\hfil\null\nobreak\hfill\copy\QEDbox\medbreak}
\newdimen\HIindent
\newbox\HIbox
\def\setHI#1{\setbox\HIbox=\hbox{#1}\HIindent=\wd\HIbox}
\def\HI#1{\par\hangindent\HIindent\hangafter=0\noindent\leavevmode
\llap{\hbox to\HIindent{#1\hfil}}\ignorespaces}

\baselineskip12pt
\parskip2.5pt plus 1pt
\hyphenation{Hei-del-berg}
\def\L|Abk:#1|Sig:#2|Au:#3|Tit:#4|Zs:#5|Bd:#6|S:#7|J:#8||{%
\edef\TEST{[#2]}
\expandafter\ifx\csname#1\endcsname\TEST\relax\else
\immediate\write16{#1 hat sich geaendert!}\fi
\expandwrite\AUX{\neverexpand\ref{#1}{\TEST}}
\HI{[#2]}
\ifx-#3\relax\else{#3}: \fi
\ifx-#4\relax\else{#4}{\sfcode`.=3000.} \fi
\ifx-#5\relax\else{\it #5\/} \fi
\ifx-#6\relax\else{\bf #6} \fi
\ifx-#8\relax\else({#8})\fi
\ifx-#7\relax\else, {#7}\fi\Par}

\def\B|Abk:#1|Sig:#2|Au:#3|Tit:#4|Reihe:#5|Verlag:#6|Ort:#7|J:#8||{%
\edef\TEST{[#2]}
\expandafter\ifx\csname#1\endcsname\TEST\relax\else
\immediate\write16{#1 hat sich geaendert!}\fi
\expandwrite\AUX{\neverexpand\ref{#1}{\TEST}}
\HI{[#2]}
\ifx-#3\relax\else{#3}: \fi
\ifx-#4\relax\else{#4}{\sfcode`.=3000.} \fi
\ifx-#5\relax\else{(#5)} \fi
\ifx-#7\relax\else{#7:} \fi
\ifx-#6\relax\else{#6}\fi
\ifx-#8\relax\else{ #8}\fi\Par}

\def\Pr|Abk:#1|Sig:#2|Au:#3|Artikel:#4|Titel:#5|Hgr:#6|Reihe:{%
\edef\TEST{[#2]}
\expandafter\ifx\csname#1\endcsname\TEST\relax\else
\immediate\write16{#1 hat sich geaendert!}\fi
\expandwrite\AUX{\neverexpand\ref{#1}{\TEST}}
\HI{[#2]}
\ifx-#3\relax\else{#3}: \fi
\ifx-#4\relax\else{#4}{\sfcode`.=3000.} \fi
\ifx-#5\relax\else{In: \it #5}. \fi
\ifx-#6\relax\else{(#6)} \fi\PrII}
\def\PrII#1|Bd:#2|Verlag:#3|Ort:#4|S:#5|J:#6||{%
\ifx-#1\relax\else{#1} \fi
\ifx-#2\relax\else{\bf #2}, \fi
\ifx-#4\relax\else{#4:} \fi
\ifx-#3\relax\else{#3} \fi
\ifx-#6\relax\else{#6}\fi
\ifx-#5\relax\else{, #5}\fi\Par}
\setHI{[KKLV]\ }
\sfcode`.=1000

\def\rho{\varrho}

\def\alp{\alpha}
\def\lam{\lambda}
\def\ebar{{\overline\eta}}
\def\E{{E_\eta}}

\def\F{{F_\eta}}
\def\Fa{{F_\eta^{(\alpha)}}}
\def\Ja{{J_\lambda^{(\alpha)}}}

\def\P{{P_\lambda}}
\def\J{{J_\lambda}}
\def\Oa{{\Omega}}
\fontdef Ss:cmss10,,.
\font\BF=cmbx10 scaled \magstep1
\font\CSC=cmcsc10 

\baselineskip15pt

{\baselineskip1.5\baselineskip\rightskip0pt plus 5truecm
\leavevmode\vskip0truecm\noindent
\BF A new scalar product for nonsymmetric Jack polynomials
}
\vskip1truecm
\leftline{{\CSC Siddhartha Sahi}%
\footnote*{\rm This work was supported by an NSF grant,
and undertaken at IAS, Princeton. }}
\leftline{Department of Mathematics, Rutgers University, 
New Brunswick NJ 08903, USA
}
\bigskip
\beginsection 1. Introduction

Jack polynomials are a remarkable family of polynomials in 
$n$ variables $x=(x_1,\cdots,x_n)$  with coefficients in the 
field $\FF:=\QQ(\alp)$ where $\alp$ is an indeterminate. 
They arise naturally in several statistical, 
physical, combinatorial, and representation theoretic 
considerations. 

The symmetric polynomials (\cite{M1}, \cite{St}, 
\cite{LV}, \cite{KS}) $\J=\Ja$ are indexed by partitions 
$\lam=(\lam_1,\cdots,\lam_n)$ where $\lam_1\ge\cdots\ge\lam_n\ge0$. 
The nonsymmetric polynomials $\F=\Fa$ (\cite{Op}, \cite{KS}, 
and \S2) are indexed by compositions $\eta=(\eta_1,\cdots,\eta_n)$ 
where $\eta_i\ge0$ are integers. They constitute orthogonal bases, 
respectively, for symmetric polynomials and all polynomials, with 
respect the scalar product $\<f,g\>_0\equiv \int_T f(x) \overline{g(x)}
d_\alp x$ where $d_\alp x$ is the measure
$\prod_{i\ne j}(1-x_ix_j^{-1})^{1/\alp}dx$ 
on the $n$-torus $T\equiv\{x\in\CC^n\mid |x_i|=1\}$.  


Now the $\J$ are also orthogonal for a second (``combinatorial'')
scalar product defined
by setting $\<m_\lam,g_\mu\>_s=\delta_{\lam\mu}$, where $m_\lam(x)$
is the symmetrized monomial obtained by summing $x^\eta:= 
x_1^{\eta_1}\cdots x_n^{\eta_n}$ over all distinct permutations
$\eta$ of $\lam$; and $g_\lam$ is defined as the coefficient
of $m_\lam(y)$ in the expansion
$$ \prod_{i,j=1}^n{1\over(1-x_i y_j)^{1/\alp}}
=\sum_\lam g_\lam(x)m_\lam(y).$$  

The two scalar products, although related via the Dyson conjecture,
(\cite{M1} VI.9), are really quite different in nature. The first 
has been generalized by Macdonald to arbitrary root systems. The 
second has intimate connections with the combinatorics of the symmetric 
group but, to our knowledge, no analogs for other Weyl groups. 

Following Dunkl \cite{Du}, we define a nonsymmetric analog  
of the second product by setting $\<x^\eta ,q_\gamma\>= 
\delta_{\eta\gamma}$, where the polynomials $q_\eta$ are 
defined by the expansion
$$
\prod_{i=1}^n{1\over(1-x_i y_i)}
\prod_{i,j=1}^n{1\over(1-x_i y_j)^{1/\alp}}
=\sum_\eta q_\eta(x)y^\eta.$$  

\Theorem dunkl. The polynomials $F_\eta$ are mutually 
orthogonal with respect to $\<,\>$.

We also derive explicit formulas for the norm $\<F_\eta,F_\eta\>$
and the function value $F_\eta(1^n)$. 

To motivate our results, we review the corresponding
formulas for $J_\lam$: Recall that
the {\it diagram\/} of a partition $\lam$ 
consists  of the points $(i,j)$ in $\ZZ^2$ with $1\le i\le n$
and $1\le j\le\lam_i$, drawn in ``matrix fashion'' with $i$ 
increasing from top to bottom and $j$ from left to right.
For $s=(i,j)\in \lam$, the number of points to the right, 
left, below, and above $s$ are called its
{\it arm\/} $a(s)\equiv\lam_i-j$, {\it coarm\/} 
$a'(s)\equiv j-1$, {\it leg\/} 
$l(s)\equiv\#\{k>i\mid j\le \lam_k\}$ and 
{\it coleg\/} $l'(s)\equiv i-1$.

Define $b_\lam=\prod_s b(s)$, $c_\lam=\prod_s c(s)$, 
and $c'_\lam=\prod_s c'(s)$,  
where $b(s)= \alp a'(s)+ n-l'(s)$, 
$c(s)=\alp a(s)+ l(s)+1$ and $c'(s)=\alp(a(s)+1)+ l(s)$. Then Stanley has 
shown \cite{St}:
$$j_\lam\equiv\<J_\lam,J_\lam\>_s =c_\lam c'_\lam 
\quad\|and|\; 
J_\lam(1^n)= b_\lam.$$

For a composition $\eta$, we define its diagram, arms, and coarms 
just as for a partition. However the leg and 
coleg at $s=(i,j)\in\eta$ are defined to be sums of their
``lower'' and ``upper'' parts. Thus we write 
$l(s)= ul(s)+ll(s)$ and 
$l'(s)=ul'(s)+ll'(s)$ where
$$ \eqalign{
ll(s)=\#\{k>i \mid j\le \eta_k\le\eta_i\};
\quad ul(s)&=\#\{k<i \mid j\le \eta_k+1\le\eta_i\}; \cr
ll'(s)=\#\{k>i \mid \eta_k>\eta_i\};\quad  
ul'(s)&=\#\{k<i \mid \eta_k\ge\eta_i\}.}$$
(If $\eta$ {\it is\/} a partition then $l(s)$ and $l'(s)$ 
agree with the previous definition.)

We now define constants  
$d_\eta=\prod_{s\in\eta}d(s)$,
$d'_\eta=\prod_{s\in\eta}d'(s)$, 
and $e_\eta=\prod_{s\in\eta}e(s)$, 
where
$d(s)  =\alp (a(s)+ 1)+l(s)+1$;
$d'(s) =\alp (a(s)+1)+l(s)$; 
and $e(s) =\alp (a'(s)+1)+n-l'(s)$.

With this notation, our second main result is:

\Theorem coef. $f_\eta\equiv\<F_\eta,F_\eta\>=d_\eta d'_\eta$.

In the course of proving this, we obtain short, new proofs 
of the following fundamental results which, since they are not 
directly related to the combinatorial product $\< , \>$, can also 
be derived from more general results of Macdonald \cite{M2} and 
Cherednik \cite{C2}:

\Theorem spec. $\F(1^n)= e_\eta$.

\Theorem symm. $j_\lam^{-1}\J =
\sum_\eta f_\eta^{-1}\F$,
summed over all distinct permutations $\eta$ of $\lam$.

Here is a brief sketch of the argument. 
First of all, the $\E$ are simultaneous eigenfunctions for 
a family of differential-reflection operators --- the Cherednik 
operators (see \cite{C1}, \cite{KS}, \cite{Op}), and we 
deduce \cite{dunkl} from properties of these operators. 
Next, we use the action of the symmetric group (as described in 
\cite{KS}) to derive recursion lemmas for the constants
$e_\eta$, $d_\eta$, and $d'_\eta$.
This proves \cite{spec} and reduces \cite{coef} 
to the case where $\eta$ is a partition. We complete the
proof by deducing Theorems \ncite{coef} and \ncite{symm} 
from Stanley's results in the symmetric case.

\beginsection 2. Preliminaries

For the reader's convenience we provide here an abbreviated 
introduction to some of the known results on Jack polynomials, 
emphasizing those needed in subsequent sections.

The symmetric Jack polynomial $\J$ is characterizable by two properties: 
first, it is an eigenfunction of a certain (Debiard--Sekiguchi)
differential operator; second, it is 
of the form $c_\lam m_\lam+$ ``lower terms in the dominance order''. 
Here  ``lower'' means an $\FF$-combination of $m_\mu$ 
such that $\mu_1+\cdots+\mu_i\le \lam_1+\cdots+\lam_i$ 
for all $i$, with equality for $i=n$. 

The nonsymmetric Jack polynomial $\F$ is similarly characterizable: 
first, it is an eigenfunction 
of the {\it Cherednik\/} operators $\xi_1,\cdots,\xi_n$, and second it
is of the form $d_\eta x^\eta+$ ``lower terms in the dominance order''.

The dominance order on compositions (and monomials) is defined 
by combining the order on partitions with the Bruhat order on the 
symmetric group $S_n$ (all written $\ge$): 
For each composition $\eta$ there is a unique {\it minimal\/} 
permutation $w_\eta$ such $\eta^+\equiv w_\eta^{-1}\eta$ is a partition
(also unique). 
We say that $\eta\ge \zeta$ if either $\eta^+>\zeta^+$, or if 
$\eta^+=\zeta^+$ but $w_\eta\le w_\zeta$. 

The $\xi_i$ are defined as follows: Let $s_{ij}$
be the operator which interchanges $x_i$ and $x_j$.
Then $N_{ij}= {1\over x_i - x_j}(1-s_{ij})$ preserves 
$\FF[x]$, and  $\xi_i=\xi^x_i:=
\alpha x_i\partial_i+ \sum_{j<i}N_{ij}x_j+\sum_{j>i}x_jN_{ij}.$
If we write $\rho=(0,-1,\cdots,-n+1)$, and put $\ebar=
\alp\eta+w_{\eta}\rho$, then $\F$ can also be characterized
up to a multiple by the eigen-equations $\xi_i\F= \ebar_i\F$. 
As noted in \cite{KS}, the eigenvalues can also be written as 
$\ebar_i=\alp\eta_i -(k'_i+k''_i)$ where
$$ k'_i=\#\{k<i\mid \eta_k\ge\eta_i\};\quad
k''_i=\#\{k>i\mid \eta_k>\eta_i\}.$$

The polynomials $\J$ and $\F$ are ``integral'' forms of the 
Jack polynomials. However we shall also find it convenient
to work with the ``monic'' forms $\P=\J/c_\lam$ and $\E=
\F/d_\eta$.

The following proposition, (Cor. 4.2, Lemma 2.4, and 
Prop. 4.3 in \cite{KS}), describes 
the action on $\E$ of the transpositions $s_i\equiv s_{i\,i+1}$, 
and of the operator $\Phi$ which acts on compositions by $\Phi(\eta):= 
(\eta_2,\cdots,\eta_n,\eta_1+1)$, and on functions by 
$\Phi f(x_1,\cdots,x_n)=x_n f(x_n,x_1,\cdots,x_{n-1})$.

\Proposition prel.
\item{1)} $E_{\Phi\eta} =\Phi\E$.
\item{2)} if $\eta_i=\eta_{i+1}$ then $s_i\E=E_{s_i\eta}.$ 
\item{3)} if $\eta_i>\eta_{i+1}$ then $d\E=(ds_i+1)E_{s_i\eta}$ 
with $d=\ebar_i-\ebar_{i+1}$.

We shall also need the following weak form of \cite{symm} (see
\cite{KS}, Lemma 2.3):

\Proposition prel2.  $\P$ is a linear combination
of $\E$ such that $\eta^+=\lam.$

\beginsection 3. The Cauchy formula

By a simple argument (\cite{St}, Prop. 2.1), the assertion 
$\<\J,J_\mu\>_s=j_\lam\delta_{\lam\mu}$ is equivalent to the 
``Cauchy'' formula $ \prod_{i,j=1}^n{1\over(1-x_i y_j)^{1/\alp}}
=\sum_\lam j_\lam^{-1} \J(x)\J(y).$  
By the same argument, to establish $\<\F,F_\gamma\>=
f_\eta\delta_{\eta\gamma}$ we need 
to prove the nonsymmetric Cauchy formula 
$\Oa=\sum_\eta f_\eta^{-1} F_\eta(x) F_\eta(y)$ for
$$\Oa\equiv \prod_{i=1}^n{1\over(1-x_i y_i)}
\prod_{i,j=1}^n{1\over(1-x_i y_j)^{1/\alp}}.$$

We start by showing:

\Lemma. For each $i=1,\cdots,n$ we have $\xi^x_i\Oa=\xi^y_i\Oa$.

\Proof: Let $\tau$ be the operator on $\FF(x,y)$ which interchanges
$x_i$ and $y_i$ for each $i$. Since $\tau(\Oa)=\Oa$,
and $\tau(\xi^x_i\Oa)=\xi^y_i\Oa$, 
it suffices to prove that ${\xi^x_i\Oa\over \Oa}$ is $\tau$-invariant. 

Since $N_{ij}x_j=x_iN_{ij}-1$ we get
$\xi^x_i=\alpha x_i\partial_i+
\sum_{j<i}x_iN_{ij}+\sum_{j>i}x_jN_{ij}-(i-1)$.

$$
\|Now|\quad 
{\alp x_i\partial_i\Oa\over \Oa}=
 -{\alp x_i y_i\over 1-x_iy_i } - \sum_j {x_iy_j\over 1-x_iy_j}= 
 -{\alp x_i y_i\over 1-x_iy_i } +n- \sum_j {1\over 1-x_iy_j},$$ 
$$
\|and|\quad {N_{ij}\Oa\over\Oa}= {1\over x_i-x_j} 
\left[1-{(1-x_iy_i)(1-x_jy_j) \over(1-x_iy_j)(1-x_jy_i)}
\right]={y_i-y_j\over (1-x_iy_j)(1-x_jy_i)}. $$  

So, writing $f \equiv g$ to denote equivalence modulo 
$\tau$-invariant terms, we get 
$$\eqalign{
{\xi^x_i\Oa\over \Oa}&\equiv 
-\sum_{j\ne i} {1\over 1-x_iy_j}
+\sum_{j<i}{x_i(y_i-y_j)\over (1-x_iy_j)(1-x_jy_i)}  
+\sum_{j>i}{x_j(y_i-y_j)\over (1-x_iy_j)(1-x_jy_i)}\cr   
&\equiv 
-\sum_{j\ne i} {1\over 1-x_iy_j}
-\sum_{j<i}{1\over1-x_jy_i}  
+\sum_{j>i}{1\over1-x_iy_j}\cr
&\equiv 
-\sum_{j<i}\left( {1\over1-x_iy_j} +{1\over1-x_jy_i} \right)\equiv0. 
}$$

The lemma follows. \qed

\Proof: (of \cite{dunkl}). As noted above, it suffices to prove the 
existence of an expansion of the form
$\Oa=\sum_\eta f_\eta^{-1} F_\eta(x) F_\eta(y)$ for
some constants $f_\eta$.

Expanding $\Oa$ in terms of the monomial
basis $x^\eta$, and re-expressing in terms of $\F(x)$
we deduce that there are polynomials $C_\eta$ in $\FF(y)$ such that 
$\Oa=\sum_\eta\F(x)C_\eta(y).$

Applying the lemma to this expression, and comparing coefficients of $\F(x)$ 
we obtain the equations $\xi^y_iC_\eta(y)=\ebar_iC_\eta(y)$.
It follows that $C_\eta$ is a multiple of $\F$. \qed

\beginsection 4. The recursion lemmas

In this section we determine how $e_\eta$, $d_\eta$ and $d'_\eta$ 
transform under $s_i$ and $\Phi$

\Lemma Philem.  For all compositions $\eta$, we have 
$d_{\Phi\eta}/d_\eta =\ebar_1+\alp+n= e_{\Phi\eta}/e_\eta$.

\Proof: We think of the diagram of $\Phi\eta$ as obtained 
from that of $\eta$ by removing the top row, moving all the other
rows up one unit, adding a point to the top row and
appending this row to the bottom of the diagram.

To prove the first equality, we consider the new point to 
have been {\it prefixed\/} to the top row of $\eta$.
It is easy to check from the definitions that the arms and legs 
of all other points are unchanged after this 
procedure, while for the new point $s=(n,1)\in\Phi\eta$ 
we have $a(s)= 
\eta_1$ and $l(s)= \#\{k>1\mid \eta_k \le \eta_1\}$. 
This gives ${d_{\Phi\eta}/ d_\eta}= 
 \alp(a(s)+1)+ l(s)+1=\ebar_1+\alp+n.$ 

For the second equality, we consider the new point to have 
been {\it appended\/} to the top row of $\eta$. This time the
coarms and colegs of other points stay the same, and for the new point
$s=(n,\eta_1+1)\in\Phi\eta$ 
we have $a'(s)= \eta_1$ and $l'(s)= 
\#\{k>1\mid \eta_k > \eta_1\}$. Thus
${e_{\Phi\eta}/ e_\eta}= \alp(a'(s)+1)+
n-l'(s)=\ebar_1+\alp+n.  $  \qed

\Lemma silem. We have $e_{s_i\eta}=e_\eta$ for all $\eta$. Also
if $\eta_i>\eta_{i+1}$ with  $d=\ebar_i-\ebar_{i+1}$,
then ${d_{s_i\eta}/ d_\eta}= {(d+1)/ d}$ and
 ${d'_{s_i\eta}/ d'_\eta}= {d/(d-1)}.$

\Proof: Comparing definitions of $l'(s)$ and 
$\ebar$ we see that $e_\eta\equiv \prod_s[\alp 
(a'(s)+1)+n-l'(s)]$ can be rewritten
as $\prod_{i=1}^n\prod_{j=1}^{\eta_i}(n+\ebar_i-j\alp )$.
Since $\ebar=w_\eta(\alp\eta^++\rho)$ is a permutation of
$\overline{s_i\eta}=w_{s_i\eta}(\alp\eta^++\rho)$,
from the last expression it follows that $e_\eta=e_{s_i\eta}$.

For the assertions about $d_\eta$ and $d'_\eta$, note that
$d=\alp(\eta_i-\eta_{i+1}) + 
(k_{i+1}'-k_i')+(k_{i+1}''-k_i'').$ 
Moreover, for the point $s=(i,\eta_{i+1}+1)\in\eta$ we have
$\eta_i-\eta_{i+1}=a(s)+1$ and
$$\eqalign{ k_{i+1}'-k'_i
&= \#\{k<i+1\mid \eta_k\ge\eta_{i+1}\}-
    \#\{k<i\mid \eta_k\ge\eta_i\}\cr
&=1+\#\{k<i\mid \eta_i>\eta_k\ge\eta_{i+1}\}=1+ul(s).}
$$
Similarly $k_{i+1}''-k''_i=
ll(s)$, and we deduce that $d=\alp(a(s)+1)+ l(s)+1$.

Now the diagram of $s_i\eta$ is obtained from 
that of $\eta$ by switching rows $i$ and $i+1$. This increases 
the leg of the point $s=(i,\eta_{i+1}+1)\in\eta$ by $1$, while 
other legs and all arms remain the same. Thus
$$ d_{s_i\eta}/d_\eta = {\alp(a(s)+1)+ (l(s)+1)+1\over
 \alp(a(s)+1)+ (l(s)+1)}= {d+1\over d}.$$
The assertion about $d'_\eta$ follows similarly. \qed

\beginsection 5. The explicit formulas

We first establish \cite{spec}.

\Proof: (of \cite{spec}) Specializing $x=1^n$ in \cite{prel} 
we get the recursions

\item{1)} $E_{\Phi\eta}(1^n)=E_\eta(1^n)$ for all $\eta$, and
\item{2)} if $\eta_i>\eta_{i+1}$ then
$E_{s_i\eta}(1^n) ={d\over d+1}\E(1^n)$ 
with $d=\ebar_i-\ebar_{i+1}$.

By the previous lemmas, $e_\eta/d_\eta$ satisfies the same 
recursions. Checking the trivial case $\eta=(0,\cdots,0)$,
we deduce $E_\eta(1^n)=e_\eta/d_\eta$ for all $\eta$. 
The result follows.  \qed

\Lemma red. If \cite{coef} holds for partitions then it holds 
for all compositions.

\Proof: First note that the nonsymmetric Cauchy formula of \S3 
may be rewritten as $\Oa=\sum_\eta r_\eta\E(x)\E(y)$ where 
$r_\eta=f_\eta^{-1} d_\eta^2$. Next, observe that $s_i^x\Oa = 
s_i^y\Oa$, where the superscripts indicate action in the 
respective variables. Combining these, we get
$$(*)\quad\sum r_\eta (s_i\E)(x) \E(y)=\sum r_{s_i\eta}
 E_{s_i\eta}(x) (s_iE_{s_i\eta})(y).$$ 

If $\eta_i>\eta_{i+1}$ and $d=\ebar_i-\ebar_{i+1}$,  
then $ s_iE_{s_i\eta}= 
-{1\over d} E_{s_i\eta}+E_\eta$ and 
$s_i\E= {1\over d}\E+ {d^2-1\over d^2}E_{s_i\eta}$.
The first follows directly from \cite{prel} and
the second after multiplication by $ds_i-1$. 
Computing coefficients of $E_{s_i\eta}(x)\E(y)$ on both sides 
of $(*)$ we get $r_{s_i\eta} = {d^2-1\over d^2}r_\eta.$

From \cite{silem} we have $d_{s_i\eta}= {d+1\over d}d_\eta$ and
$d'_{s_i\eta}= {d\over d-1} d'_\eta$. Combining these we deduce
that $f_\eta^{-1}d_\eta d_\eta'=r_\eta d'_\eta d_\eta^{-1}$ is
$s_i$-invariant, and the result follows. \qed

To continue we need the following formula for
symmetric Jack polynomials.

\Lemma Las. We have $\displaystyle
\prod_{i=1}^n{1\over (1-x_i)^r}=\sum_\lam k_\lam {\J(x)\over j_\lam }$
with $\displaystyle k_\lam= \prod_{s\in\lam} [\alp(r+a'(s))-l'(s)].$

\Proof: Expanding the left side in terms of $m_\lam$ 
and then in terms of $\J$, we deduce that the $k_\lam$ are 
{\it polynomials\/} in $\FF[r]$. Next, by the symmetric 
Cauchy formula in $m\ge n$ variables:
$$\prod_{i,j=1}^m{1\over (1-x_iy_j)^{1/\alp}}=
\sum_{\mu} j_\mu^{-1} J_\mu(x_1,\cdots,x_m)
 J_\mu(y_1,\cdots,y_m).$$

Let $\lam=(\mu_1,\cdots,\mu_n)$ be the first $n$ components of
$\mu$. By the ``stability'' of Jack polynomials (\cite{St} 
Prop. 2.5), $J_\mu(x_1,\cdots,x_n,0,\cdots,0)$ equals
$J_\lam(x_1,\cdots,x_n)$ if $\mu_{n+1}=\cdots=\mu_m=0$
and equals {\it zero\/} otherwise.
Thus setting $x_{n+1}=\cdots=x_m=0$ and $y_1=\cdots=y_m=1$ in 
the $m$-variable formula, and using the expressions for 
$j_\mu$ and $J_\mu(1^m)$, we obtain
$$\prod_{i=1}^n{1\over (1-x_i)^{m/\alp}}=\sum_\lam 
 \prod_{s\in\lam} [\alp a'(s)+m-l'(s)] 
 j_\lam^{-1}\J(x).$$

This implies the formula for $r=m/\alp$ and, by polynomiality
of $k_\lam$, in general. \qed

Finally we prove Theorems \ncite{coef} and \ncite{symm}.

\Proof: (of Theorems \ncite{coef} and \ncite{symm}) 
Setting $y_1=\cdots=y_n=1$ in the nonsymmetric Cauchy 
formula, and using \cite{Las} and \cite{spec} to
evaluate the two sides we get
$$
\sum_\eta f_\eta^{-1} e_\eta \F=
\sum_\lam \prod_{s\in\lam}
[\alp({n\over\alp}+1+a'(s))-l'(s)] { \J \over j_\lam}
=\sum_\lam e_\lam j_\lam^{-1}\J
$$ 

Using \cite{prel2} and \cite{silem} ($e_\eta=e_{\eta^+}$),
this gives $j_\lam^{-1} J_\lam=\sum_{\eta^+=\lam} f_\eta^{-1} 
\F$ which is \cite{symm}. 
Comparing coefficients of $x^\lam$ on both sides
we get $j_\lam^{-1}c_\lam= f_\lam^{-1} d_\lam$,
which can be rewritten as $f_\lam= c'_\lam d_\lam$. Since
$c'_\lam=d'_\lam$, we deduce \cite{coef} for partitions 
and, by \cite{red}, in general. \qed

\beginsection References. References

\baselineskip12pt
\parskip2.5pt plus 1pt
\hyphenation{Hei-del-berg}
\def\L|Abk:#1|Sig:#2|Au:#3|Tit:#4|Zs:#5|Bd:#6|S:#7|J:#8||{%
\edef\TEST{[#2]}
\expandafter\ifx\csname#1\endcsname\TEST\relax\else
\immediate\write16{#1 hat sich geaendert!}\fi
\expandwrite\AUX{\neverexpand\ref{#1}{\TEST}}
\HI{[#2]}
\ifx-#3\relax\else{#3}: \fi
\ifx-#4\relax\else{#4}{\sfcode`.=3000.} \fi
\ifx-#5\relax\else{\it #5\/} \fi
\ifx-#6\relax\else{\bf #6} \fi
\ifx-#8\relax\else({#8})\fi
\ifx-#7\relax\else, {#7}\fi\Par}

\def\B|Abk:#1|Sig:#2|Au:#3|Tit:#4|Reihe:#5|Verlag:#6|Ort:#7|J:#8||{%
\edef\TEST{[#2]}
\expandafter\ifx\csname#1\endcsname\TEST\relax\else
\immediate\write16{#1 hat sich geaendert!}\fi
\expandwrite\AUX{\neverexpand\ref{#1}{\TEST}}
\HI{[#2]}
\ifx-#3\relax\else{#3}: \fi
\ifx-#4\relax\else{#4}{\sfcode`.=3000.} \fi
\ifx-#5\relax\else{(#5)} \fi
\ifx-#7\relax\else{#7:} \fi
\ifx-#6\relax\else{#6}\fi
\ifx-#8\relax\else{ #8}\fi\Par}

\def\Pr|Abk:#1|Sig:#2|Au:#3|Artikel:#4|Titel:#5|Hgr:#6|Reihe:{%
\edef\TEST{[#2]}
\expandafter\ifx\csname#1\endcsname\TEST\relax\else
\immediate\write16{#1 hat sich geaendert!}\fi
\expandwrite\AUX{\neverexpand\ref{#1}{\TEST}}
\HI{[#2]}
\ifx-#3\relax\else{#3}: \fi
\ifx-#4\relax\else{#4}{\sfcode`.=3000.} \fi
\ifx-#5\relax\else{In: \it #5}. \fi
\ifx-#6\relax\else{(#6)} \fi\PrII}
\def\PrII#1|Bd:#2|Verlag:#3|Ort:#4|S:#5|J:#6||{%
\ifx-#1\relax\else{#1} \fi
\ifx-#2\relax\else{\bf #2}, \fi
\ifx-#4\relax\else{#4:} \fi
\ifx-#3\relax\else{#3} \fi
\ifx-#6\relax\else{#6}\fi
\ifx-#5\relax\else{, #5}\fi\Par}
\setHI{[ABC]\ }
\sfcode`.=1000

\L|Abk:C1|Sig:C1|Au:Cherednik, I.|Tit:A unification of the
Knizhnik\_Zamolodchikov equations and Dunkl operators via
affine Hecke algebras|Zs:Invent. Math.|Bd:106|S:411--432|J:1991||

\L|Abk:C2|Sig:C2|Au:Cherednik, I.|Tit:Nonsymmetric Macdonald
polynomials|Zs:IMRN|Bd:10|S:483--515|J:1995||

\L|Abk:Du|Sig:Du|Au:Dunkl, C.|Tit:Intertwining operators and polynomials 
associated with the symmetric group|Zs:-|Bd:-|S:preprint|J:1996||

\L|Abk:KS|Sig:KS|Au:Knop, F; Sahi, S.|Tit:A recursion and a combinatorial 
formula for Jack polynomials|Zs:Invent. Math.|Bd:-|S:-|J:to appear||

\L|Abk:LV|Sig:LV|Au:Lapointe, L.; Vinet, L.|Tit:Exact operator
solution of the Calogero-Sutherland model|Zs:-|Bd:-|S:-|J:to appear||

\B|Abk:M1|Sig:M1|Au:Macdonald, I.|Tit:Symmetric Functions and Hall
Polynomials|Reihe:2nd ed.|Verlag:Clarendon Press|Ort:Oxford|J:1995||

\L|Abk:M2|Sig:M2|Au:Macdonald, I.|Tit:Affine Hecke algebras
and orthogonal polynomials|Zs:S\'em. Bourbaki|Bd:797|S:1--18|J:1995||

\L|Abk:Op|Sig:Op|Au:Opdam, E.|Tit:Harmonic analysis for certain
representations of graded Hecke algebras%
|Zs:Acta Math.|Bd:175|S:75--121|J:1995||

\L|Abk:St|Sig:St|Au:Stanley, R.|Tit:Some combinatorial properties of
Jack symmetric functions|Zs:Adv. Math.|Bd:77|S:76--115|J:1989||

\bye